\documentclass{amsart}
\usepackage{amsmath}
\usepackage{amssymb}
\usepackage{a4}		
\title{The distance of a permutation from a subgroup of $S_n$}
\author{Richard G.E. Pinch}
\address{2 Eldon Road, Cheltenham, Glos GL52 6TU, U.K.}
\email{rgep@chalcedon.demon.co.uk}
\dedicatory{For Bela Bollobas on his 60th birthday}
\date{20 November 2005}
\newtheorem{theorem}{Theorem}
\newtheorem{proposition}[theorem]{Proposition}
\newtheorem{problem}{Problem}

\newcommand{\del}{\partial}

\newcommand{\suchthat}{\ \vert\ }

\def\abs|#1|{{\left\vert{#1}\right\vert}}
\def\O(#1){\mbox{O}{\left({#1}\right)}}
\def\C(#1){C{\left({#1}\right)}}
\def\GF(#1){{\mathbb{F}_{#1}}}			
\def\gen<#1>{{\left\langle{#1}\right\rangle}}
\def\paren(#1){{\left({#1}\right)}}

\begin{document}

\begin{abstract}
We show that the problem of computing the distance of a given permutation from 
a subgroup $H$ of $S_n$ is in general NP-complete, even under the restriction 
that $H$ is elementary Abelian of exponent 2.  The problem is shown to be 
polynomial-time equivalent to a problem related to finding a maximal partition 
of the edges of an Eulerian directed graph into cycles and this problem is in
turn equivalent to the standard NP-complete problem of Boolean satisfiability.
\end{abstract}

\maketitle

\suppressfloats

\section{Introduction}

We show that the problem of computing the distance of a given permutation 
from a subgroup $H$ of $S_n$ is in general NP-complete, even under the 
restriction that $H$ is elementary Abelian of exponent 2.  The problem is 
polynomial-time equivalent to finding a
maximal partition of the edges of an Eulerian directed graph into cycles and
this is turn equivalent to the standard NP-complete problem {\bf 3-SAT}.

\section{   Distance in the symmetric group }

We define {\it Cayley distance} in a symmetric group as the minimum number of
transpositions which are needed to change one permutation to another by
post-multiplication
$$
d(\rho,\pi) = \min \left\lbrace n \suchthat \rho\tau_1\ldots\tau_n = \pi,
                   \quad\hbox{$\tau_i$ transpositions}\ \right\rbrace .
$$
It is well-known that Cayley distance is a metric on $S_n$ and that
it is homogeneous, that is, 
$d(\rho,\pi) = d\left(I,\rho^{-1}\pi\right)$.
Further, the distance of a permutation $\pi$
from the identity in $S_n$ is $n$ minus the number of cycles in $\pi$.

If $H$ is a subgroup of $S_n$, then we define the distance of a permutation
$\pi$ from $H$ as
$$
d(H,\pi) = \min_{\eta \in H} d(\eta,\pi) .
$$

We refer to Critchlow \cite{Cri:partrank} and Diaconis \cite{Dia:grpreps} for background and
further material on the uses of the Cayley and other metrics on $S_n$.

\begin{problem}[Subgroup--Distance]
\par\noindent{\sc  Instance:} Symmetric group $S_n$, element $\pi \in S_n$, 
elements $\left\lbrace h_1 , \ldots , h_r\right\rbrace$ of $S_n$, integer $K$.
\par\noindent{\sc  Question:} Is there an 
element $\eta \in H = \langle h_1 , \ldots , h_r \rangle$ such 
that $d(\eta,\pi) \le K$? 
\end{problem}

The natural measure of this problem is $nr$ where $r$ is the length of the
list of generators.  The following result shows that every subgroup of $S_n$ 
has a set of generators of length at most $n^2$ and hence we are
justified in taking $n$ as the measure of the
various problems derived from {\bf Subgroup--Distance}.

\begin{proposition}
Every subgroup of $S_n$ can be generated by at most $n^2$ elements.
\end{proposition}
\begin{proof}
Let $H$ be a subgroup of $S_n$.  It is clear that $H$ is generated by the 
union of one Sylow subgroup for every prime $p$ dividing the order $\#H$.
The order of a Sylow $p$-subgroup of $H$ is $p^b$ where $p^b$ divides $\#H$ and
hence $n!$.  It is well-known (e.g.~Dickson \cite{Dic:history} I, chap 9)
that the power of $p$ dividing $n!$ is at 
most $\frac{n}{p-1}$, so $b \le n$.  But consideration of the composition factors
shows that a $p$-group of order $p^b$ can be generated by a set of at most $b$ 
elements: hence any Sylow subgroup of $H$ can be generated by at most
$n$ elements.  Further, the set of prime factors of the order of $H$
forms a subset of the set of prime factors of $n!$, that is, of the primes up
to $n$, and there are at most $n$ such primes.
Hence $H$ can be generated by a set of at most $n^2$ elements.

\end{proof}

Although we do not need the stronger result, it can be shown that any
subgroup of $S_n$ can be generated by at most $3n-2$ elements. 

We define a subset of $S_n$ to 
be {\it involutions with disjoint support} 
{\it (IDS)} to be set of elements of the
form $\gamma_j = \left(x_j^{(1)}\,y_j^{(1)}\right) \ldots 
                     \left(x_j^{(r_j)}\,y_j^{(r_j)}\right)$
where the $x_j^{(i)}, y_j^{(i)}$ are all distinct.
The subgroup generated by an IDS is clearly elementary Abelian with exponent 2.
Define the {\it width} of an IDS to be the maximum number of 2-cycles $r_j$
in the generators $\gamma_j$.  
The problem {\bf IDS$w$--Subgroup--Distance} is the 
problem {\bf Subgroup--Distance}
with the list of generators restricted to be an IDS of width at most $w$.

\begin{theorem}
The problem {\bf IDS$6$--Subgroup--Distance} is NP-complete.
\end{theorem}

The Theorem will follow from combining Theorem 3 and Theorem 7.
We deduce immediately that the more general
problem {\bf Subgroup--Distance} is also NP-complete.

By contrast, the problem of deciding whether the distance is zero, that is, 
testing for membership of a subgroup of $S_n$, has a polynomial-time
solution, an algorithm first given by Sims \cite{Sims:compmeth} and shown to have a
polynomial-time variant by Furst, Hopcroft and Luks \cite{FHL:polyperm}.  See  Babai,
Luks and Seress \cite{BLS:fastperm} and Kantor and Luks \cite{KL:compquot} for a survey of related 
results.

\section{   Switching circuits }

Let $G=(V,E)$ be a directed graph with vertex set $V$ and edge set $E$.  
(We allow loops and multiple edges.)
For each vertex $v$ define $e_+(v)$ to be the set of edges out of $v$
and $e_-(v)$ the set of edges into $v$.  The in-valency $\del_-(v) = \#e_-(v)$
and the out-valency $\del_+(v) = \#e_+(v)$.
We define a {\it switching circuit} to be a directed graph $G$ for which
$\del_+(v) = \del_-(v) = \del(v)$, say, and for which 
there is a labelling $l_\pm(v)$ of each set $e_\pm(v)$ with the integers from
$1$ to $\del(v)$.  (The labels at each end of an edge are not related.)
A {\it routing} $\rho$ for a switching circuit is a 
choice of permutation $\rho(v) \in S_{\del(v)}$ for each vertex $v$.  
Clearly there is a correspondence between routings for a switching circuit $G$
and decompositions of the edge set of $G$ into directed cycles.
We define a {\it polarisation} $T$ for a switching circuit $G$ to be an 
equivalence relation on the set of vertices such that equivalent vertices 
have the same valency, and call $(G,T)$ a {\it polarised switching circuit}.
We say that a routing $\rho$ {\it respects} the polarisation $T$ 
if the permutations $\rho(x)$ and $\rho(y)$ are equal
whenever $x$ and $y$ are equivalent vertices under $T$.
We shall sometimes refer to a switching circuit without a polarisation,
or with a polarisation for which all the classes are trivial, 
as {\it unpolarised}.

\begin{problem}[Polarised--Switching--Circuit--Maximal--Routing]
\par\noindent{\sc  Instance:} Polarised switching circuit $(G,T)$,
positive integer $K$.
\par\noindent{\sc  Question:} Is there a routing which respects $T$ and has
at least $K$ cycles in the associated edge-set decomposition?
\end{problem}

We define the {\it width} of a polarisation to be the maximum 
number of vertices in an equivalence class of $T$.  
The problem {\bf Width$w$--Valency$v$--Maximal--Routing} is the 
problem {\bf Polarised--Switching--Circuit--Maximal--Routing}
with the width of $T$ constrained to be
at most $w$ and the in-{} and out-valency of each vertex in $V$ 
constrained to be at most $v$.

\begin{theorem}
Problem {\bf Width$6$--Valency$2$--Maximal--Routing} is NP-complete.
\end{theorem}

\section{    Proof of Theorem 3 }

We shall show that the problem {\bf 3-SAT}, [LO2] of Garey and 
Johnson \cite{GJ:NP}, which is known to be NP-complete, can be reduced to
the problem {\bf Width$6$--Valency$2$--Maximal--Routing}.

We define a polarised switching circuit $(G,T)$
to be {\it Boolean} (or {\it binary})
if every vertex has in- and out-valency 1 or 2.  
To each class $C$ of the
polarisation $T$ we associate a Boolean variable $a(C)$.  There is then
a 1-1 correspondence between routings $\rho$ which respect $T$ and 
assignments of
truth values to the variables $a(C), C \in T$ by specifying that $a(C)$
is 0 (false) if and only if the permutation $\rho(v)$ is the identity in
$S_2$ for every $v$ in $C$, and 1 (true) if and only if $\rho(v) = (1\,2)$.

We denote a vertex in a polarisation class associated with the Boolean 
variable $a$ as in Figure 1.  Out convention for drawing the diagrams
will be to assume the edges round each vertex labelled so that 1 is
denoted by either ``straight through'' or ``turn right''.

\setlength{\unitlength}{0.5pt}
\begin{figure}[ht]
\begin{picture}(100,100)(-50,-50)
\put(-5,-5){$a$}
\put(0,0){\circle{60}}
\put(-50,50){\vector(1,-1){35}}
\put(-50,-50){\vector(1,1){35}}
\put(15,15){\vector(1,1){35}}
\put(15,-15){\vector(1,-1){35}}
\put(-50,50){1}
\put(50,-50){1}
\put(-50,-50){2}
\put(50,50){2}
\end{picture}
\ \ \ %
\begin{picture}(100,100)(-50,-50)
\put(-50,50){\vector(1,-1){100}}
\put(-50,-50){\line(1,1){35}}
\put(15,15){\vector(1,1){35}}
\end{picture}
\ \ \ %
\begin{picture}(100,100)(-50,-50)
\qbezier(-50,50)(0,0)(50,50)
\put(49,49){\vector(1,1){1}}
\qbezier(-50,-50)(0,0)(50,-50)
\put(49,-49){\vector(1,-1){1}}
\end{picture}

\vspace{0.2in}

\begin{picture}(100,100)(-50,-50)
\put(-5,-5){$a$}
\put(0,0){\circle{60}}
\put(-15,15){\vector(-1,1){35}}
\put(-50,-50){\vector(1,1){35}}
\put(50,50){\vector(-1,-1){35}}
\put(15,-15){\vector(1,-1){35}}
\put(-50,50){2}
\put(50,-50){1}
\put(-50,-50){1}
\put(50,50){2}
\end{picture}
\ \ \ %
\begin{picture}(100,100)(-50,-50)
\qbezier(-50,50)(0,0)(50,50)
\put(-49,49){\vector(-1,1){1}}
\qbezier(-50,-50)(0,0)(50,-50)
\put(49,-49){\vector(1,-1){1}}
\end{picture}
\ \ \ %
\begin{picture}(100,100)(-50,-50)
\qbezier(-50,-50)(0,0)(-50,50)
\put(-49,49){\vector(-1,1){1}}
\qbezier(50,50)(0,0)(50,-50)
\put(49,-49){\vector(1,-1){1}}
\end{picture}

\caption{A vertex in a switching circuit
associated with the Boolean variable $a$,
and the routings with $a = 1$ and $a = 0$ respectively}

\end{figure}
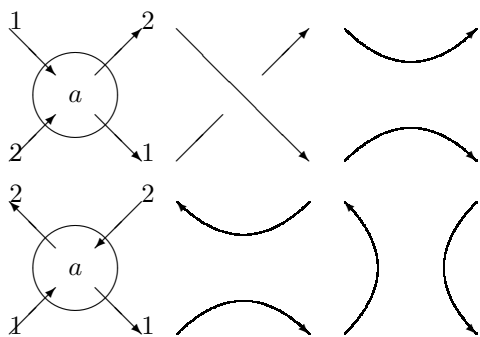

We associate a vertex with the negated variable $\bar a$ by exchanging the
input labels 1 and 2.

Our proof will proceed by finding polarised switching circuits for which the 
number of maximal cycles in a routing is a Boolean function of the variables.

For a single Boolean variable $a$ define $I(a)$ to be the switching circuit
in Figure 2.

\begin{figure}[ht]
\begin{picture}(160,160)(-80,-80)
\put(-5,-5){$a$}
\put(0,0){\circle{60}}
\put(-15,15){\vector(-1,1){35}}
\put(-50,-50){\vector(1,1){35}}
\put(50,50){\vector(-1,-1){35}}
\put(15,-15){\vector(1,-1){35}}
\qbezier(-50,50)(-60,60)(-50,70)
\qbezier(-50,-50)(-60,-60)(-50,-70)
\qbezier(50,50)(60,60)(50,70)
\qbezier(50,-50)(60,-60)(50,-70)
\qbezier(-50,70)(0,80)(50,70)
\qbezier(-50,-70)(0,-80)(50,-70)
\end{picture}
\caption{The switching circuit $I(a)$.}

\end{figure}
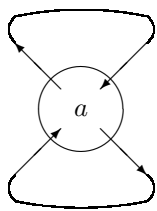

For a pair of Boolean variables $(a,b)$ define the polarised
switching circuit $E(a,b)$ as in Figure 3.

\begin{figure}[ht]
\begin{picture}(240,160)(-80,-80)
\put(-5,-5){$a$}
\put(0,0){\circle{60}}
\put(-15,15){\vector(-1,1){35}}
\put(-50,-50){\vector(1,1){35}}
\put(50,50){\vector(-1,-1){35}}
\put(15,-15){\vector(1,-1){35}}
\qbezier(-50,50)(-60,60)(-50,70)
\qbezier(-50,-50)(-60,-60)(-50,-70)
\qbezier(-50,70)(50,80)(150,70)
\qbezier(-50,-70)(50,-80)(150,-70)
\put(95,-5){$b$}
\put(100,0){\circle{60}}
\put(85,15){\vector(-1,1){35}}
\put(50,-50){\vector(1,1){35}}
\put(150,50){\vector(-1,-1){35}}
\put(115,-15){\vector(1,-1){35}}
\qbezier(150,50)(160,60)(150,70)
\qbezier(150,-50)(160,-60)(150,-70)
\end{picture}
\caption{The switching circuit $E(a,b)$.}
\end{figure}
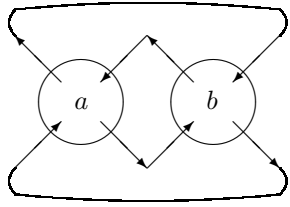

Further define the polarised switching circuit $F(a,b)$ as in Figure 4.

\begin{figure}[hb]
\begin{picture}(240,240)(-80,-180)
\put(-105,-5){$a$}
\put(-100,0){\circle{60}}
\put(-150,50){\vector(1,-1){35}}
\put(-150,-50){\vector(1,1){35}}
\put(-85,15){\vector(1,1){35}}
\put(-85,-15){\vector(1,-1){35}}
\qbezier(-150,50)(-160,60)(-150,70)
\qbezier(-150,70)(0,80)(150,70)
\put(95,-5){$b$}
\put(100,0){\circle{60}}
\put(50,50){\vector(1,-1){35}}
\put(50,-50){\vector(1,1){35}}
\put(115,15){\vector(1,1){35}}
\put(115,-15){\vector(1,-1){35}}
\qbezier(150,50)(160,60)(150,70)
\put(-5,-105){$a$}
\put(0,-100){\circle{60}}
\put(-50,-50){\vector(1,-1){35}}
\put(-50,-150){\vector(1,1){35}}
\put(15,-85){\vector(1,1){35}}
\put(15,-115){\vector(1,-1){35}}
\qbezier(-50,-150)(-60,-160)(-50,-170)
\qbezier(-50,-170)(100,-180)(250,-170)
\put(195,-105){$b$}
\put(200,-100){\circle{60}}
\put(150,-50){\vector(1,-1){35}}
\put(150,-150){\vector(1,1){35}}
\put(250,-50){\vector(-1,-1){35}}
\put(215,-115){\vector(1,-1){35}}
\qbezier(250,-150)(260,-160)(250,-170)
\put(-50,50){\vector(1,0){100}}
\put(50,-150){\vector(1,0){100}}
\put(-150,-50){1}
\put(250,-50){1}
\end{picture}
\caption{The switching circuit $F(a,b)$.}
\end{figure}
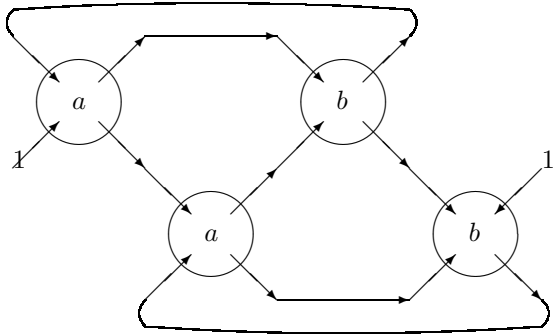

Define $G(a,b)$ to be the disjoint union of $F(a,b)$ and $E(\bar a,b)$.

\begin{proposition}
\begin{enumerate}
\item The number of cycles in a routing for $I(a)$ is 2 if $a=1$ and 
otherwise 1.
\item The number of cycles in a routing for $E(a,b)$ is 2 if $a=b$ and
otherwise 1.
\item The number of cycles in a routing for $F(a,b)$ is 2 if $a \ne b$, 
3 if $a=b=1$ and 1 if $a=b=0$.
\item The number of cycles in a routing for $G(a,b)$ is 2 if $a=b=0$ and
4 otherwise.
\end{enumerate}
\end{proposition}
 
\begin{proof}
In each case we simply enumerate the cases.

\end{proof}
 
For a triple of Boolean variables $(a,b,c)$ define the polarised
switching circuit $A(a,b,c)$ as in Figure 5.

\setlength{\unitlength}{0.3pt}
\begin{figure}[ht]
\begin{picture}(800,400)(0,-400)
\put(0,0){
\begin{picture}(600,100)(-150,50)
\put(-105,-5){$a$}
\put(-100,0){\circle{60}}
\put(-150,50){\vector(1,-1){35}}
\put(-150,-50){\vector(1,1){35}}
\put(-85,15){\vector(1,1){35}}
\put(-85,-15){\vector(1,-1){35}}
\put(95,-5){$b$}
\put(100,0){\circle{60}}
\put(50,50){\vector(1,-1){35}}
\put(50,-50){\vector(1,1){35}}
\put(115,15){\vector(1,1){35}}
\put(115,-15){\vector(1,-1){35}}
\put(295,-5){$c$}
\put(300,0){\circle{60}}
\put(250,50){\vector(1,-1){35}}
\put(250,-50){\vector(1,1){35}}
\put(315,15){\vector(1,1){35}}
\put(315,-15){\vector(1,-1){35}}
\put(-50,50){\line(1,0){100}}
\put(150,50){\line(1,0){100}}
\put(-150,50){1}
\put(-150,-50){2}
\put(350,50){4}
\end{picture}}
\put(100,-100){
\begin{picture}(600,100)(-150,50)
\put(-105,-5){$a$}
\put(-100,0){\circle{60}}
\put(-150,50){\vector(1,-1){35}}
\put(-150,-50){\vector(1,1){35}}
\put(-85,15){\vector(1,1){35}}
\put(-85,-15){\vector(1,-1){35}}
\put(95,-5){$b$}
\put(100,0){\circle{60}}
\put(50,50){\vector(1,-1){35}}
\put(50,-50){\vector(1,1){35}}
\put(115,15){\vector(1,1){35}}
\put(115,-15){\vector(1,-1){35}}
\put(295,-5){$c$}
\put(300,0){\circle{60}}
\put(250,50){\vector(1,-1){35}}
\put(250,-50){\vector(1,1){35}}
\put(315,15){\vector(1,1){35}}
\put(315,-15){\vector(1,-1){35}}
\put(-150,-50){3}
\put(350,50){3}
\end{picture}}
\put(200,-200){
\begin{picture}(600,100)(-150,50)
\put(-105,-5){$a$}
\put(-100,0){\circle{60}}
\put(-150,50){\vector(1,-1){35}}
\put(-150,-50){\vector(1,1){35}}
\put(-85,15){\vector(1,1){35}}
\put(-85,-15){\vector(1,-1){35}}
\put(95,-5){$b$}
\put(100,0){\circle{60}}
\put(50,50){\vector(1,-1){35}}
\put(50,-50){\vector(1,1){35}}
\put(115,15){\vector(1,1){35}}
\put(115,-15){\vector(1,-1){35}}
\put(295,-5){$c$}
\put(300,0){\circle{60}}
\put(250,50){\vector(1,-1){35}}
\put(250,-50){\vector(1,1){35}}
\put(315,15){\vector(1,1){35}}
\put(315,-15){\vector(1,-1){35}}
\put(-150,-50){4}
\put(350,50){1}
\end{picture}}
\put(300,-300){
\begin{picture}(600,100)(-150,50)
\put(-105,-5){$a$}
\put(-100,0){\circle{60}}
\put(-150,50){\vector(1,-1){35}}
\put(-150,-50){\vector(1,1){35}}
\put(-85,15){\vector(1,1){35}}
\put(-85,-15){\vector(1,-1){35}}
\put(95,-5){$b$}
\put(100,0){\circle{60}}
\put(50,50){\vector(1,-1){35}}
\put(50,-50){\vector(1,1){35}}
\put(115,15){\vector(1,1){35}}
\put(115,-15){\vector(1,-1){35}}
\put(295,-5){$c$}
\put(300,0){\circle{60}}
\put(250,50){\vector(1,-1){35}}
\put(250,-50){\vector(1,1){35}}
\put(315,15){\vector(1,1){35}}
\put(315,-15){\vector(1,-1){35}}
\put(-50,-50){\line(1,0){100}}
\put(150,-50){\line(1,0){100}}
\put(-150,-50){5}
\put(350,50){5}
\put(350,-50){2}
\end{picture}}

\end{picture}

\caption{The switching circuit $A(a,b,c)$.}
\end{figure}
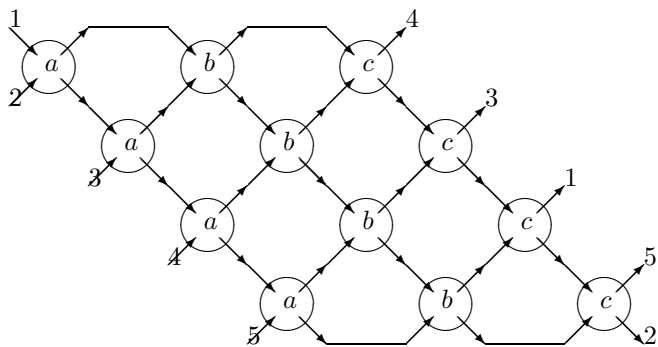


\begin{proposition}
The number of cycles in a routing for $A(a,b,c)$ is 1 if $a=b=c=0$ and
3 otherwise.
\end{proposition}

\begin{proof}
Again, in each case we simply enumerate the cases.

\end{proof}

\begin{theorem}
There is a polynomial-time parsimonious transformation from the problem
{\bf 3-SAT} to the problem {\bf Width$6$--Valency$2$--Maximal--Routing}.
\end{theorem}
\begin{proof}
Suppose we have an instance of {\bf 3-SAT}: that is, a Boolean formula $\Phi$
of length $l$ in variables $x_i$ which is a conjunct of $k$ clauses each of 
which is a disjunct of at most three variables (possibly negated).
We transform $\Phi$ into a formula $\Phi'$ in variables $y_i^j$ 
by replacing the $j^{\text{th}}$ occurence of 
variable $x_i$ by the variable $y_i^j$ and conjoining clauses
$\left(y_i^1 \equiv y_i^2\right) \wedge \ldots \wedge
 \left(y_i^{(r_{i-1})} \equiv y_i^{(r_i)}\right)$
where the variable $x_i$ occurs $r_i$ times in $\Phi$.
Clearly $\Phi$ and $\Phi'$ represent the same Boolean function and have the
same number of satisfying assignments.  Every
variable in $\Phi'$ occurs at most three times, and at most once in a disjunct
deriving from a clause in $\Phi$.  Let $n$ be the total number of variables
in $\Phi'$; certainly $n \le l$.

We form a polarised switching circuit $\Psi$ from $\Phi'$ as follows.
Take a circuit $B(x,y,z)$ for every clause
in $\Phi'$ of the form $(x \vee y \vee z)$; take a circuit $G(x,y)$
for every clause 
in $\Phi'$ of the form $(x \vee y)$; take a circuit $I(x)$ for every clause
in $\Phi'$ of the form $(x)$; take a circuit $E(x,y)$ for every clause 
in $\Phi'$ of the form $(x \equiv y)$.  Let the number of circuits of types
$B$, $G$, $I$ and $E$ taken to form $\Psi$ be $b$, $g$, $i$, 
and $e$ respectively.  Put $M = 3b + 4g + 2i + 2e$.
The resulting polarised switching circuit has $n$ classes, and each class in
the polarisation is involved in at most one circuit of the 
form $B$, $G$ or $I$: hence each class contains at most $4 + 1 + 1 = 6$ 
vertices and the number of vertices in $\Psi$ is thus at most $6n$.
Furthermore, a routing for $\Psi$ has
$M$ cycles if and only if the corresponding assignment of Boolean values 
gives $\Phi'$, and hence $\Phi$, the value 1; otherwise a routing has less 
than $M$ cycles.

\end{proof}

Since the problem {\bf 3-SAT} is known to be NP-complete,
we immediately deduce that the problem 
{\bf Width$16$--Valency$2$--Maximal--Routing} is NP-complete as well.
This proves Theorem 3.

\section{    Switching circuits and IDS }

In this section we obtain a polynomial-time equivalence between the problems
{\bf Width$w$--Valency$2$--Maximal--Routing} and 
{\bf IDS$w$--Subgroup--Distance}.

\begin{theorem}
There is a polynomial-time parsimonious equivalence between problems
{\bf Width$w$--Valency$2$--Maximal--Routing} and 
{\bf IDS$w$--Subgroup--Distance}.
\end{theorem}
\begin{proof}
Suppose we have an instance of {\bf IDS$w$--Subgroup--Distance},
that is, an element $\pi$ of $S_n$ together with an IDS on a set of $t$ 
generators $\left\lbrace\gamma_j\right\rbrace$ with
$\gamma_j = \left(x_j^{(1)}\,y_j^{(1)}\right) \ldots 
               \left(x_j^{(r_j)}\,y_j^{(r_j)}\right),$ 
the $x_j^{(i)}, y_j^{(i)}$ all distinct and all the $r_j \le w$.  
We construct a polarised switching circuit on a graph, vertex 
set $V = \left\lbrace P(1), \ldots, P(n)\right\rbrace \cup 
             \left\lbrace Q(1,1), \ldots, Q(t,r_t)\right\rbrace$.
Each vertex $P(k)$ will be of in-valence and out-valence 1;
each vertex $Q(j,i)$ will be of in-valence and out-valence 2.
For each $j$ up to $t$ and $i$ up to $r_j$ we take edges
from $P(x_j^i)$ and from $P(y_j^i)$ to $Q(j,i)$ labelled 1 and 2 respectively,
and edges from $Q(j,i)$ to $P\left(\pi\left(x_j^i\right)\right)$ 
and to $P\left(\pi\left(y_j^i\right)\right)$ again labelled
1 and 2 respectively.  We define a polarisation $T$ on $V$ by taking $t$ 
classes $C_j = \left\lbrace Q(j,i) \suchthat i = 1, \ldots, r_j\right\rbrace$;
clearly the width of $T$ is at most $w$.

Conversely, suppose
we have an instance of {\bf Width$w$--Valency$2$--Maximal--Routing},
that is, a directed graph $(V,E)$ with every vertex $v$ having in-{} and 
out-valency two, a 
labelling $l_\pm(v) : e_\pm(v) \rightarrow \{ 1,2 \}$
of edges into and out of each vertex $v$, and an equivalence 
relation $T$ on $V$ with $t$ classes each of size at most $w$.
Put $n = \#E$.  We define a permutation $\pi$ of $E$ as follows.
For an edge $e$ into a vertex $v$, let $\pi(e)$ be the edge $f$ out of $v$
which has label $l_{+}(v)(f)$ equal to $l_{-}(v)(e)$.  We further define an
IDS by writing down a set of 
generators $\left\lbrace\gamma_j\right\rbrace$ as follows.  For each class of 
vertices $C_j = \left\lbrace v_i^j \suchthat i = 1,\ldots, r_j \right\rbrace$ 
in the polarisation $T$, let $\gamma_j$ be the product of transpositions of 
the form $\left(f_i^j g_i^j\right)$ where $f_i^j$ and $g_i^j$ are the edges 
out of vertex $v_i^j$. Since each class in $T$ has at most $w$ elements, each 
generator $\gamma_j$ is composed of at most $w$ transpositions.

In each case there is a correspondence between routings $\rho$ of the
switching circuit which respect the polarisation $T$
and permutations of the form $\pi\eta$ where $\eta$
runs over the elements of the subgroup $H$ of $S_n$ generated by
the $\gamma_j$: in this correspondence the number of cycles in the 
routing $\rho$ is equal to the number of cycles in the permutation $\pi\eta$.
Hence $\pi$ is within
distance $d$ of the group generated by the $\gamma_j$ if and only if there is 
a routing $\rho$ with at least $n-d$ cycles.

\end{proof}




\providecommand{\bysame}{\leavevmode\hbox to3em{\hrulefill}\thinspace}
\providecommand{\MR}{\relax\ifhmode\unskip\space\fi MR }
\providecommand{\MRhref}[2]{%
  \href{http://www.ams.org/mathscinet-getitem?mr=#1}{#2}
}
\providecommand{\href}[2]{#2}

\end{document}